# A Non-Pythagorean Musical Scale Based on Logarithms

Robert P. Schneider
Department of Mathematics • University of Kentucky
Lexington, KY 40506, USA
`robert.schneider@uky.edu`

## Abstract

A new musical scale devised by the author, based on natural logarithms, is described. Most of the logarithmic pitches bear no correspondence to the twelve tones of the ancient tuning system attributed to Pythagoras, based on ratios of whole numbers, nor to the chromatic tones of scales in equal temperament used widely in the modern era. Logarithms obey a special arithmetic compared to whole and rational numbers, which can be heard in beat frequencies between tones of the scale, suggesting extended harmonic possibilities by incorporating the beat frequencies into compositions. The author uses the broad term "non-Pythagorean" to describe the logarithmic musical scale, as the ratios of pitches are usually irrational numbers.

## Introduction

An interesting musical scale may be defined by letting successive tones have the ratios to one another of the natural logarithms of successive whole numbers. While melodies in the scale sound somewhat exotic, novel musical possibilities arise from mathematical properties of logarithms, particularly with respect to the construction of chords.

We use the broad term "non-Pythagorean" to describe this scale, as the logarithmic pitches bear no correspondence to the tuning system based on ratios of whole numbers attributed to the Greek mathematician Pythagoras [1], nor to the chromatic tones of scales in equal temperament [3], but rather are based on a sequence of irrational intervals. Non-Pythagorean scales are used in traditional music around the world, such as in Indian classical music and in gamelan orchestras of Indonesia; however, the scale described here follows more closely in the experimental footsteps of twentieth century composers such as Ivan Wyschnegradsky, Karlheinz Stockhausen, Iannis Xenakis and La Monte Young, who have utilized mathematically-based scales in their compositions [5].

Successive pitches in the logarithmic musical scale grow closer and closer together, and the number of distinct tones in each octave increases exponentially with each successive octave. When listening to the lower octaves of the scale, the ear strains to associate the tones with traditional chromatic pitches lying near them. However, the sequence of pitches does not correspond to the traditional twelve tones. As successive tones grow closer together they become nearly indistinguishable, superficially resembling microtonal tunings, and the resemblance to a traditional scale disappears altogether.

## Mathematical definition of the scale

To begin the logarithmic musical scale on middle C we define a fundamental pitch $f$ by

$$f = \frac{c}{\log 4} \approx 190.44 \text{ Hz}$$

where $c = 264$ Hz is a typical frequency used for middle C, the notation $\log M$ denotes the natural logarithm of a whole number $M$, and we use the approximation $\log 4 \approx 1.3863$. We note that any frequency may be substituted for 264 Hz as the value of $c$ in this expression, to produce a sequence of tones starting on a different pitch from middle C. We then define and notate the $M$th pitch in the logarithmic scale as

$$M^* = f \log M.$$

We will adopt the convention of indicating a logarithmic tone with an asterisk adjoined to the number of the tone in the sequence of pitches, e.g. 10* denotes the 10th tone in the sequence. If we let $M = 1$ in the above formula, we find the first tone 1* is effectively silent, as $f \log 1 = 0$. The second tone 2* corresponds to C one octave below middle C on the piano, the fourth tone 4* is middle C itself, and the sixteenth tone 16* corresponds to C one octave above middle C.

The eighth tone 8* lies a perfect fifth above tone 4*, corresponding to the pitch G in the chromatic scale, as $8^*/4^* = \log 8 / \log 4 = 3/2$. More generally, tone $X^*$ is a perfect fifth above tone $Y^*$ whenever $X = n^3$ and $Y = n^2$ for some integer $n$, as $X^*/Y^* = 3/2$. However, for $M < 16$ every other ratio between pitches is an irrational number; the tones do not lie in the chromatic scale or correspond to any pitches traditionally played on a piano keyboard.

As noted above, successive pitches grow closer together as $M$ increases due to the slow growth rate of the logarithm function, and each octave contains a greater number of tones than the preceding octave. The C two octaves above middle C is tone 64*, thus there are 48 tones in the octave beginning with tone 16*. The next C is tone 256*, thus there are 192 tones in the octave beginning with 64*. The fourth octave contains 768 tones, and generally the $n$th octave of the scale contains $3 \cdot 2^{2n}$ tones.

The construction of a piano with such a great number of keys in a single octave is impractical, if not impossible. The electronic medium seems best suited to composing with the logarithmic scale, if use of an extended range of tones is desired; however, a restricted set of tones can be played on certain traditional instruments, as noted below.

### A twelve-tone scale suitable for playing on a traditional keyboard

To generate a twelve-tone subset of the scale suitable for playing upon a re-tuned traditional piano keyboard, and suitable for composing using traditional notation, we associate tone $4^* = f \log 4$ with both the pitch and keyboard position middle C on the keyboard; and we associate subsequent logarithmic pitches with the subsequent tones in the octave, such that the formula for the $N$th tone in the octave above middle C is given by

$$f \log(3 + N),$$

where $N$ ranges from 1 to 12. For the present purpose it is a useful coincidence that there exist twelve integers ranging from 4 to 15, as well as that $\log 16 = 2 \log 4$, for a twelve-tone octave of the present scale may be associated with the integers on this interval, and tone 16* lies precisely one octave above tone 4*.

Below we calculate the frequencies comprising the twelve-tone scale, approximated to two decimal places. Higher and lower octaves of the resulting sequence of tones may be associated with other octaves of the keyboard by multiplying or dividing the following frequency values by powers of 2:

$$C^* = \text{tone } 4^* = f \log 4 = 264 \text{ Hz}$$
$$C\#^* = \text{tone } 5^* = f \log 5 = 306.49 \text{ Hz}$$
$$D^* = \text{tone } 6^* = f \log 6 = 341.22 \text{ Hz}$$
$$D\#^* = \text{tone } 7^* = f \log 7 = 370.57 \text{ Hz}$$
$$E^* = \text{tone } 8^* = f \log 8 = 396 \text{ Hz}$$
$$F^* = \text{tone } 9^* = f \log 9 = 418.43 \text{ Hz}$$
$$F\#^* = \text{tone } 10^* = f \log 10 = 438.49 \text{ Hz}$$
$$G^* = \text{tone } 11^* = f \log 11 = 456.64 \text{ Hz}$$

$$G\#* = \text{tone } 12* = f \log 12 \doteq 473.22 \text{ Hz}$$
$$A* = \text{tone } 13* = f \log 13 \doteq 488.46 \text{ Hz}$$
$$A\#* = \text{tone } 14* = f \log 14 \doteq 502.57 \text{ Hz}$$
$$B* = \text{tone } 15* = f \log 15 \doteq 515.71 \text{ Hz}.$$

In the above list we employ the convention of adjoining an asterisk to the traditional name of the keyboard position of a tone, to indicate a logarithmic note. Then the logarithmic notes C*, D*, E*, F#*, etc. fall on the keys traditionally associated with those letters.

It would be desirable to play an acoustic piano tuned to the pitches calculated above, such as the Fluid Piano built by British inventor Geoff Smith. This author has composed a number of short pieces [8, 9] in this twelve-tone scale using a sine wave generator [7], as well as a rock song incorporating the scale [10], released on albums by his pop group The Apples in Stereo; and has also composed pieces using the H-π Tuning Box for MIDI synthesizers designed by Aaron Andrew Hunt, using a re-tuned kalimba, or African thumb piano, and using a fretless box guitar with the positions corresponding to the above pitches marked off on the fingerboard.

### Difference tones and ideas for future compositions

A beat frequency is a frequency defined as the absolute difference between two pitches sounded simultaneously; when this frequency is in the audible range and is perceived as a distinct pitch, it is called a difference tone. While audible beat frequencies are often considered undesirable artifacts in music, due to the subtraction property $\log x - \log y = \log(x/y)$ of logarithms there is potential in the logarithmic scale to compose with not only the tones themselves, but with the beat frequencies within chords.

Specifically, if whole-numbered tones $X*$ and $Y*$ are chosen such that $X/Y$ is an integer, the beat frequency created by sounding the two pitches together will itself be a tone $(X/Y)*$ of the logarithmic scale, allowing an additional degree of manipulation of tones within the scale—that is, purposeful melodic lines may be constructed in the realm of difference tones.

The author wishes to explore such avenues for composition more deeply, having worked out only rudimentary music-theoretic and number-theoretic observations at the present time. One interesting question that presents itself is whether optimal sets of integers may be chosen such that pairs of corresponding logarithmic pitches will most often produce difference tones that are also members of the resulting scale. Another question is whether second-order beat frequencies between difference tones are themselves detectible by the ear, as well as third-order beat frequencies and so on, for these might also be used in compositions. There are subtleties of the mathematics of beat frequencies that we have neglected to detail here—for instance, the relative amplitudes and phases of the waveforms come into play, as well as physical properties of the resonators [6]—but one may assume ideal conditions when composing.

The author also intends to produce longer works exploiting the special properties of the scale, such as a string quartet currently in progress, to be performed by violins, viola and cello with the finger positions corresponding to logarithmic pitches marked on the fingerboards. Specially-constructed flutes, xylophones and marimbas, specially-tuned harps, and other acoustic instruments might also be utilized in performances. Certainly the natural harmonics of each instrument will interfere to some degree with the effects of the beat frequencies, but such considerations need not limit the composer's imagination. Specifics regarding the constructions and tunings of such instruments will be detailed in a future, more thorough treatise; the requisite principles of instrument construction—e.g. the locations of holes in flutes or of frets on stringed instruments—are well-known and easily extended for the present purpose [4].

We should note that a similar scale may be defined using the square root function instead of the logarithm function, which also contains a twelve-tone octave playable on a traditional keyboard; however, the square root function does not yield special properties with respect to the beat frequencies produced.

## Concluding remarks

Melodies and harmonies obtained from the logarithmic scale have unusual qualities by comparison with those produced in the chromatic scale, hinting at new avenues for musical expression. For example, logarithmic chords produced by a sine wave generator possess distinctive textures due to beat frequencies, at times resembling crickets chirping, the ringing of bells and other environmental sounds. Listeners have noted a melancholy, "ethereal" [2] quality to the compositions produced thus far.

There is an alien beauty to this non-Pythagorean musical scale, stemming from the mathematics of logarithms, once the listener becomes accustomed to the strange intervals.

## References

[1] J. Murray Barbour, *Tuning and Temperament: A Historical Survey*, Mineola, New York: Dover Publications, Inc. (2004).

[2] Mick Hamer, *Flexible Scales and Immutable Octaves*, New Scientist, Issue 2644 (Feb. 23, 2008), pp. 32–34.

[3] Leon Harkleroad, *The Math Behind the Music*, New York: Cambridge University Press (2006).

[4] Bart Hopkin, *Musical Instrument Design*, Tuscon, Arizona: See Sharp Press (1996).

[5] Andy Mackay, *Electronic Music*, Minneapolis, Minnesota: Control Data Publishing (1981).

[6] Haynes R. Miller, *18.03 Differential Equations Course Notes, Spring 2010*, Massachusetts Institute of Technology: MIT OpenCourseWare, http://ocw.mit.edu/courses/mathematics/18-03-differential-equations-spring-2010/readings/supp_notes (May 23, 2012), Web.

[7] Stefan Scheer, *The Sound of Tomorrow: Robert Schneider und die Apples in Stereo*, Innovation Stuntmen: Strategien fuer das 21. Jahrhundert, http://www.innovationstuntmen.com/?p=841 (Oct. 23, 2010), Web.

[8] The Apples in Stereo, "Non-Pythagorean Composition 1," *New Magnetic Wonder*, Simian/Yep Roc Records (2007), CD.

[9] The Apples in Stereo, "Non-Pythagorean Composition 3," *New Magnetic Wonder*, Simian/Yep Roc Records (2007), CD.

[10] The Apples in Stereo, "CPU," *Travellers in Space and Time*, Simian/Yep Roc Records (2010), CD.